\documentclass[11pt,a4paper]{article}
\usepackage{graphicx}
\usepackage{amsfonts}
\usepackage{latexsym}
\usepackage{epsfig}

\newcommand{\beql}[1]{\begin{equation}\label{#1}}
\newcommand{\eeq}{\end{equation}}
\newcommand{\comment}[1]{}

\newcommand{\eqref}[1]{{\rm (\ref{#1})}}

\newcommand{\Abs}[1]{{\left|{#1}\right|}}

\newcommand{\Norm}[1]{{\left\|{#1}\right\|}}

\newcommand{\Qed}{\ \\\mbox{$\Box$}}

\newcommand{\Set}[1]{{\left\{{#1}\right\}}}

\newcommand{\RR}{{\mathbb R}} 

\newcommand{\ZZ}{{\mathbb Z}}

\newcommand{\inner}[2]{{\langle #1, #2 \rangle}}

\newcommand{\dist}{{\rm dist\,}}

\newcommand{\ft}[1]{\widehat{#1}}

\newcounter{open}
\newcounter{dfn}
\def\thedfn{\arabic{dfn}}

\newcounter{obs}
\def\theobs{\arabic{obs}}

\newcounter{thm}
\newcounter{othm}
\def\theothm{\Alph{othm}}
\newenvironment{othm}{
  \sf
  \vskip 0.10in
  \refstepcounter{othm}
  \noindent{\bf Theorem\ \theothm}
}

\newcounter{mysec}
\newcounter{mysubsec}[mysec]
\newtheorem{theorem}{Theorem}
\newtheorem{corollary}{Corollary}
\newtheorem{lemma}{Lemma}

\newtheorem{remark}{Remark}

\begin{document}

\title{Distance sets corresponding to convex bodies}
\author{Mihail N. Kolountzakis}
\date{March 2003}
\maketitle

\begin{abstract}
Suppose that $K \subseteq \RR^d$ is a $0$-symmetric convex body
which defines the usual norm
$$
\Norm{x}_K = \sup\Set{t\ge 0:~~x \notin tK}
$$
on $\RR^d$.
Let also $A\subseteq\RR^d$ be a measurable set of positive upper density
$\rho$.
We show that if the body $K$ is not a polytope, or if it is a polytope with
many faces (depending on $\rho$),
then the distance set
$$
D_K(A) = \Set{\Norm{x-y}_K:~~x,y\in A}
$$
contains all points $t\ge t_0$ for some positive number $t_0$.
This was proved by Katznelson and Weiss,
by Falconer and Marstrand and by Bourgain in the case where $K$ is
the Euclidean ball in any dimension.
As corollaries we obtain (a) an extension to any dimension
of a theorem of Iosevich and \L aba
regarding distance sets with respect to convex bodies of well-distributed
sets in the plane,
and also (b) a new proof of a theorem of Iosevich, Katz and Tao about the
nonexistence of Fourier spectra for smooth convex bodies.
\end{abstract}

%%%%%%%%%%%%%%%%%%%%%%%%%%%%%%%%%%%%%%%%%%%%%
% Introduction
%%%%%%%%%%%%%%%%%%%%%%%%%%%%%%%%%%%%%%%%%%%%%
\section{Introduction}

Suppose that $K \subseteq \RR^d$ is a $0$-symmetric convex body
which defines the usual norm
$$
\Norm{x}_K = \sup\Set{t\ge 0:~~x \notin tK}
$$
on $\RR^d$.
Define the $K$-distance set of $A \subseteq \RR^d$ as the set of
$K$-distances that show up in $A$:
$$
D_K(A) = \Set{\Norm{x-y}_K:~~x,y \in A}.
$$

Katznelson and Weiss \cite{furstenberg-katznelson-weiss},
Falconer and Marstrand \cite{falconer-marstrand} (in the plane) and
Bourgain \cite{bourgain}
have proved that if $B$ is the Euclidean ball and
$A$ has positive upper density, i.e., if there is $\rho>0$ such that there are
arbitrarily large cubes $Q$ in which $A$ has a fraction at least $\rho$
of their measure: $\Abs{A\cap Q} \ge \rho\Abs{Q}$,
then $D_B(A)$ contains all numbers $t\ge t_0$, for some $t_0$.
In this paper we study the following question: 
for which other convex bodies $K$ in place of $B$
is this true?
We obtain that if this property
(of eventually all numbers showing up in $D_K(A)$)
fails for some set $A$ with positive upper density
then the body $K$ is necessarily a polytope,
and with a number of faces that is
bounded above by a number that is inversely proportional
to the density of $A$.

To state our results more precisely
let us call a measure $\mu\in M(\RR^d)$ {\em $\epsilon$-good} if
for some $R>0$ we have $\Abs{\ft{\mu}(x)} \le \epsilon$ for
$\Abs{x}\ge R$.
In what follows the Fourier Transform of a measure is defined by
$$
\ft{\mu}(\xi) = \int_{\RR^d} e^{-2\pi i \inner{x}{\xi}}~d\mu(x).
$$

\begin{theorem}\label{th:good-measure}
Suppose that $A \subseteq \RR^d$ has upper density equal to $\epsilon>0$
and that the $0$-symmetric convex body $K$ affords $C_d\epsilon$-good probability
measures supported on its boundary (the constant $C_d$ depends on the dimension only).
Then $D_K(A)$ contains all positive real numbers beyond a point.
\end{theorem}

This is complemented by the following result.
\begin{theorem}\label{th:non-polytopes-afford}
Suppose $K$ is a $0$-symmetric convex body. Then
\begin{itemize}
\item[(a)]
If $K$ is a polytope with $N$ non-parallel face directions then
it does not afford $\epsilon$-good probability measures on its boundary if
$\epsilon < 1/(\sqrt{2}N)$.
\item[(b)]
If $K$ is a polytope with $N$ non-parallel face directions then it
affords $({1\over N}+\delta)$-good probability measures on its boundary for all $\delta>0$.

If $K$ is not a polytope then it affords $\epsilon$-good probability measures
on its boundary for all $\epsilon>0$.
\end{itemize}
\end{theorem}

Theorems \ref{th:good-measure} and \ref{th:non-polytopes-afford} allow us
to prove the following result, a weaker version of which
was the motivation for this work and
was proved in dimension $2$ by Iosevich and \L aba \cite{iosevich-laba}.
A set $\Lambda\subseteq\RR^d$ is called {\em well-distributed} if there is a constant $r>0$
such that every cube of side $r$ contains at least one point of $\Lambda$.
And a set $D\in\RR$ is called separated if there exists $\epsilon>0$ such that
$\Abs{x-y}\ge\epsilon$ for any two distinct $x, y \in D$.
\begin{corollary}\label{cor:iosevich-laba-generalization}
Suppose that $\Lambda$ is a well-distributed subset of $\RR^d$ and $K$ is 
a $0$-symmetric convex body. 
If $D_K(\Lambda) \cap \RR^+$ has infinitely many gaps of length
at least $\epsilon>0$ then $K$ is a polytope.
\end{corollary}
\begin{remark}
It was proved in \cite{iosevich-laba} that if $\Lambda$ is a well-distributed
set in the plane and $D_K(\Lambda)$ is separated then $K$ is a polygon.
\end{remark}
\begin{remark}
By taking $K = (-1/2, 1/2)^d$ and $\Lambda = \ZZ^d$ we see that there
are indeed polytopes and well-distributed sets for which $D_K(\Lambda)$ is
separated. It is not clear which polytopes can play this role.
\end{remark}
\noindent
{\bf Proof of Corollary \ref{cor:iosevich-laba-generalization}.}
Suppose that the intervals $(x_k, x_k+\epsilon)$, $x_k\to+\infty$,
do not intersect $D_K(\Lambda)$.
Let $A$ be the subset of $\RR^d$ that arises if we put a copy of the body $(\epsilon/10)K$
centered at each point of $\Lambda$.
By the fact that $\Lambda$ is well-distributed we obtain that
$A$ has positive upper density and for
any two points  $x,y \in A$, we can write $x=\lambda+x_1$, $y=\mu+y_1$,
with $\lambda,\mu \in \Lambda$, $x_1, y_1 \in (\epsilon/10)K$.
It follows that $x-y = (\lambda-\mu)+(x_1-y_1)$ and
$$
{1\over 5}\epsilon \ge \Norm{x_1-y_1}_K \ge \Abs{~\Norm{x-y}_K-\Norm{\lambda-\mu}_K~},
$$
which implies that $D_K(A)$
does not intersect the intervals $(x_k+\epsilon/5, x_k+4\epsilon/5)$,
$k=1,2,\ldots$,
hence $D_K(A)$ cannot contain a half-line. From
Theorems \ref{th:good-measure} and \ref{th:non-polytopes-afford} it follows
that $K$ is a polytope, and with a number of faces
which is bounded above by a function of the upper density of $A$.
\Qed

From Corollary \ref{cor:iosevich-laba-generalization}
one can easily show that smooth convex bodies do not have Fourier spectra,
a fact first proved by Iosevich, Katz and Tao \cite{iosevich-katz-tao},
who used a different approach.
A bounded open domain $\Omega \subseteq \RR^d$ is said to have the
set $\Lambda \subseteq \RR^d$ as a {\em Fourier spectrum} if the collection
of exponentials
$$
E(\Lambda) = \Set{\exp(2\pi i \inner{\lambda}{x}):\ \lambda\in\Lambda}
$$
is an orthogonal basis for $L^2(\Omega)$.
It is easy to see from the orthogonality that any two distinct points
$\lambda$ and $\mu$ of a Fourier spectrum must satisfy
\begin{equation}\label{differences}
\ft{\chi_\Omega}(\lambda - \mu) = 0,
\end{equation}
and that the set $\Lambda$ is necessarily a well-distributed set.
\begin{corollary}
\label{cor:smooth}
{\rm (Iosevich, Katz and Tao \cite{iosevich-katz-tao})}\\
If $K$ is a smooth, $0$-symmetric convex body it does not have
a Fourier spectrum.
\end{corollary}

\noindent
{\bf Proof of Corollary \ref{cor:smooth}.}
Suppose $\Lambda$ is a Fourier spectrum of $K$.
It is a well known fact (see, for example, \cite{iosevich-katz-tao})
that if $\xi$ is a zero of $\ft{\chi_K}$ and $\xi\to\infty$
then
$$
\Norm{\xi}_{K^o} = \left({\pi\over 2} + {d\pi \over 4}\right) + k\pi + o(1),
	\ \ \ (\xi\to\infty),
$$
where $K^o$ is the dual body (which is also smooth), $d$ is the dimension
and $k$ is an integer. 

Let $R>0$ be such that any zero $\xi$ of $\ft{\chi_K}$, outside a cube of side
$R$ centered at the origin, satisfies
$$
\Norm{\xi}_{K^o} = \left({\pi\over 2} + {d\pi \over 4}\right) + k\pi + \theta,
	\ \ \ (k\in\ZZ,\ \Abs{\theta}<\pi/10).
$$
We also take $R$ to be large enough so as to be certain that we find
at least one $\Lambda$-point in any cube of side $R$.
(We can do this since $\Lambda$ is well-distributed.)

Let now the set $\Lambda'$ arise by keeping only one point of $\Lambda$
in each cube of the type $n + (-R/2, R/2)^d$, with $n\in\ZZ^d$ having
all its coordinates even. We keep nothing outside these cubes.
It follows that $\Lambda'$ is also a well
distributed set and that for any two distinct points $\lambda$ and $\mu$
of $\Lambda'$, $\mu$ is not contained in the cube of side $R$
centered at $\lambda$.
From \eqref{differences} we obtain
that for any two distinct points $\lambda,\mu\in\Lambda'$ we have
$$
\Norm{\lambda-\mu}_{K^o} = k\pi + \theta,\ \ \ (k\in\ZZ,\ \Abs{\theta}\le\pi/5).
$$
This means that the set of $K^o$ distances $D_{K^o}(\Lambda')$
has infinitely many gaps of length at least $3\pi/5$, so
by Corollary \ref{cor:iosevich-laba-generalization} $K^o$ should be a
polytope, which is a contradiction.
\Qed

%%%%%%%%%%%%%%%%%%%%%%%%%%%%%%%%%%%%%%%%%%%%%
% Proofs
%%%%%%%%%%%%%%%%%%%%%%%%%%%%%%%%%%%%%%%%%%%%%
\section{Proofs of the main theorems}
\label{sec:proofs}

\noindent
{\bf Proof of Theorem \ref{th:good-measure}.}

\noindent
{\em Notation}:
$B_1(0)$ is the unit ball in $\RR^d$ and $\omega_d$ is its volume.
In what follows the dimension $d>1$ is fixed and constants may depend on it.

As in \cite{bourgain}, whose method we follow,
it suffices to prove the following theorem.
\begin{othm}\label{th:compact}
Suppose $\epsilon, R >0$ and
let $K$ be a $0$-symmetric convex body contained in $B_1(0)$ and
$\sigma$ be a probability measure on $\partial K$ such that
\beql{compact-condition}
\Abs{\ft{\sigma}(\xi)} \le \eta(\epsilon) = {\omega_d \over 80\cdot 4^d \pi^d} \epsilon,
\ \ \ \ \mbox{if $\Abs{\xi}\ge R$}.
\eeq
Then there is $J = J(\epsilon, R)$ such that if
$A \subseteq B_1(0)$ is a set of measure $\epsilon$ and $0<t_j<1$, $j\ge 1$, is
a sequence with $t_{j+1} < {1\over 2}t_j$ then there is $j\le J$
such that
$t_j \in D_K(A)$.
\end{othm}

\noindent{\bf Theorem \ref{th:compact} implies Theorem \ref{th:good-measure}:}
Suppose that $K$ is a $0$-symmetric convex body, and $\sigma\in M(\partial K)$
is a probability measure on its surface which satisfies
$\Abs{\ft{\sigma}(\xi)} \le \eta(\epsilon)$ if $\Abs{\xi} \ge R$,
and let $A \subseteq \RR^d$ be a measurable set with upper density
larger than $\epsilon>0$.
Suppose also that Theorem \ref{th:good-measure} fails and there is
a sequence $0<x_1<x_2<\cdots$ tending to infinity all of whose elements
are not in $D_K(A)$.

Renaming, we pass to lacunary subsequence of $x_j$ such that
$x_{j+1} > 2 x_j$ with $J = J(C_d\epsilon, R)$ terms.
Let $Q$ be a cube of side-length greater than $10 x_J$ for which
$\Abs{Q \cap A} \ge \epsilon\Abs{Q}$.
Scaling everything down to $B_1(0)$ and renaming the scaled $x_j$
as $t_j$ (and reversing their order so as to have a decreasing sequence)
we now have a set of measure $\ge C \epsilon$ contained in $B_1(0)$
and a sequence $t_1>t_2>\cdots>t_J$, with $t_{j+1} < {1\over 2} t_j$.

By Theorem \ref{th:compact} there exists a $j \le J$ such that $t_j$
appears as a $K$-distance in the scaled-down set $A$, which implies
that the corresponding $x_{j'}$ appears as a $K$-distance in the set $A$,
a contradiction.
\Qed

\noindent{\bf Proof of Theorem \ref{th:compact}.}
It is enough to show that there is $j\le J(\epsilon,R)$ such that, for $t=t_j$,
$$
\int\int f(x) f(x+ty) dx d\sigma(y) > 0,
$$
where $f$ is the indicator function of $A$.
The integral on the left may be rewritten as a constant $C(t)$ times
\beql{int1}
\int\Abs{\ft{f}(\xi)}^2 \ft{\sigma}(t\xi)~d\xi.
\eeq
The integral in \eqref{int1} is broken into three parts
$$
I_1 = \int_{\Abs{\xi} \le {\delta\over t}}\cdot,\ \ \ 
I_2 = \int_{{\delta\over t} < \Abs{\xi} < {1\over \delta t}}\cdot,\ \ \ 
I_3 = \int_{{1\over \delta t}\le \Abs{\xi}}\cdot.
$$
where $\delta$ will be determined later.

We shall now need three lemmas to control the quantities $I_1, I_3$, and $I_2$.

\begin{lemma}\label{th:A1}
Let $K$ be a $0$-symmetric convex body contained in the unit ball $B_1(0)$.
Let also $\sigma$ be a probability measure on $\partial K$ and
$A \subseteq B_1(0)$ be a measurable set with indicator function $f$.
Then, writing
$$
I_1(t,\delta) = \int_{\Abs{\xi} \le {\delta\over t}}
 \Abs{\ft{f}(\xi)}^2 \ft{\sigma}(t\xi)~d\xi,
$$
we have
\beql{A1-eqn}
I_1(t,\delta) \ge {\omega_d \over 8\cdot 4^d \pi^d} \Abs{A}^2,
\ \ \ \ (\mbox{if $t\le 4 \pi \delta$}).
\eeq
\end{lemma}
{\bf Proof of Lemma \ref{th:A1}.}
Since $f$ and $\sigma$ are supported in $B_1(0)$ it follows that
$\Abs{\nabla_u\ft{f}} \le 2\pi\Abs{A}$
and
$\Abs{\nabla_u\ft{\sigma}} \le 2\pi$, for all directions $u \in S^{d-1}$.
Hence
$$
\ft{f}(\xi) \ge {\Abs{A} \over 2},\ \ \ \ \mbox{if $\Abs{\xi}\le{1\over 4\pi}$},
$$
and
$$
\ft{\sigma}(\xi) \ge {1\over 2},\ \ \ \ \mbox{if $\Abs{\xi}\le{1\over 4\pi}$}.
$$
So
\begin{eqnarray*}
I_1(t, \delta)
 & \ge & \int_{\Abs{\xi}\le{1\over4\pi}} \Abs{\ft{f}(\xi)}^2\ft{\sigma}(t\xi)~d\xi\\
 & \ge & \int_{\Abs{\xi}\le{1\over4\pi}} \left({\Abs{A}\over2}\right)^2 {1\over 2}~d\xi\\
 & \ge & {1\over 8}\Abs{A}^2 \Abs{\Set{\Abs{\xi}\le {1\over4\pi}}}\\
 & \ge & {\omega_d \over 8\cdot 4^d \pi^d} \Abs{A}^2.
\end{eqnarray*}
\Qed

\begin{lemma}\label{th:A2}
Suppose that $\sigma$ is any measure for which $\Abs{\ft{\sigma}(x)} \le \eta$
provided $\Abs{x} \ge R$ and write
$$
I_3(t, \delta) = \int_{\Abs{\xi} \ge {1 \over \delta t}}
	\Abs{g(\xi)}^2 \ft{\sigma}(t\xi)~d\xi,
$$
where $g \in L^2(\RR^d)$.
Then
\beql{A2-eqn}
\Abs{I_3(t,\delta)} \le \eta \int\Abs{g}^2,
\ \ \ \ \mbox{(if $\delta \le {1 \over R}$)}.
\eeq
\end{lemma}
{\bf Proof of Lemma \ref{th:A2}.}
\begin{eqnarray*}
\Abs{I_3(t, \delta)}
 & \le & \eta \int_{\Abs{\xi} \ge {R \over t}} \Abs{g(\xi)}^2~d\xi\\
 & \le & \eta\int\Abs{g}^2.
\end{eqnarray*}
\Qed

\begin{lemma}\label{th:A3}
If $t, \delta, \theta > 0$, $\delta<1$,
$\sigma$ is a finite measure with total variation at most $1$,
$f \in L^2(\RR^d)$ with $\epsilon = \int\Abs{f}^2$,
$$
I_2(t, \delta) = \int_{{\delta\over t} < \Abs{\xi} < {1 \over \delta t}}
			\Abs{\ft{f}(\xi)}^2\ft{\sigma}(t\xi)~d\xi,
$$
and $\Set{t_j}_{j=1}^\infty$ is a sequence with
$$
0<t_j<1 \mbox{\ and\ } t_{j+1} < {1\over 2} t_j
$$
then there is an index
$$
j \le {2 \over \log 2}~\theta^{-1} {1\over \epsilon}\log{1\over\delta}
$$
such that
\beql{A3-eqn}
\Abs{I_2(t,\delta)} \le \theta \epsilon^2.
\eeq
\end{lemma}
{\bf Proof of Lemma \ref{th:A3}.}
Note first that $\log{1\over t_{j+1}} - \log{1\over t_j} \ge \log 2$.
For $x \in \RR$ let $N(x)$ be the number of intervals
$\left({\delta\over t_j}, {1\over\delta t_j}\right)$
to which $x$ belongs.
It follows that
$$
N(x) \le {2\over\log 2}~\log{1\over\delta}
$$
as
$x \in \left({\delta\over t_j}, {1\over\delta t_j}\right)$ is equivalent to
$$
\log{1\over t_j}-\log{1\over\delta}
 < \log x < \log{1\over t_j} + \log{1\over\delta}.
$$
For any positive integer $J$ we have thus
$$
\sum_{j=1}^J \Abs{I_2(t_j, \delta)} \le
 \int_\RR \Abs{\ft{f}(\xi)}^2 \ft{\sigma}(t\xi) N(\Abs{\xi})~d\xi \le
 {2\over \log 2}~\log{1\over\delta} \cdot \epsilon,
$$
by $\int\Abs{\ft{f}(\xi)}^2~d\xi = \epsilon$ and $\Abs{\ft{\sigma}(\xi)}\le 1$.
Hence, if we let
$$
J = {2 \over \log 2}~\theta^{-1}{1\over\epsilon}\log{1\over\delta},
$$
we obtain that there is $j\le J$ for which
$$
\Abs{I_2(t_j, \delta)} \le \theta \epsilon^2.
$$
\Qed

\noindent{\bf Proof of Theorem \ref{th:compact} (continued).}
Let $\delta = 1/R$ and $t \le 4\pi\delta$.
By Lemma \ref{th:A1} we get
\beql{I1-bound}
I_1(t,\delta) \ge {\omega_d \over 8\cdot4^d\pi^d} \epsilon^2.
\eeq
By Lemma \ref{th:A2}, applied to $g = \ft f$ we get
\beql{I3-bound}
\Abs{I_3(t,\delta)} \le \eta \epsilon
	= {\omega_d \over 80\cdot 4^d\pi^d}\epsilon^2.
\eeq
Inequalities \eqref{I1-bound} and \eqref{I3-bound} hold for all $t\le 4\pi/R$.

Define $j_0 = j_0(R)$ by $t_{j_0} \le 4\pi / R$
(clearly $j_0 \le C \log R$, as $t_1\le 1$)
and $\theta = {\omega_d \over 80\cdot 4^d \pi^d}$, and
apply Lemma \ref{th:A3} to the sequence $t_{j_0}, t_{j_0+1}, \ldots$.
It follows that there is $j$ with
$$
j \le j_0 + 10\theta^{-1}{1\over \epsilon}\log{1\over\delta}
 \le C {\log R \over \epsilon} =: J(\epsilon, R)
$$
such that
\beql{I2-bound}
\Abs{I_2(t_j,\delta)} \le \theta \epsilon^2.
\eeq
Putting together \eqref{I1-bound}, \eqref{I3-bound} and \eqref{I2-bound}
we obtain for this $j$
\begin{eqnarray*}
I(t_j) &\ge& I_1(t_j,\delta) - \Abs{I_2(t_j,\delta)} - \Abs{I_3(t_j,\delta)}\\
 &\ge& {\omega_d \over 8\cdot 4^d\pi^d}\epsilon^2 - {\omega_d \over 80\cdot 4^d\pi^d}\epsilon^2
		- {\omega_d \over 80\cdot 4^d\pi^d}\epsilon^2\\
 &\ge& {\omega_d \over 40\cdot 4^d\pi^d}\epsilon^2\\
 &>& 0,
\end{eqnarray*} 
which shows that $t_j \in D_K(A)$, as we had to show.
The proof of Theorem \ref{th:compact} and therefore of Theorem \ref{th:good-measure}
is complete.
\Qed

We denote by $S^{d-1}$ the surface of the unit ball in $\RR^d$ and whenever
$\Omega$ is a hypersurface in $\RR^d$ we denote by $\sigma_\Omega$ its surface measure.
Also, if $x, y \in S^{d-1}$, by $\dist(x,y)$ we understand their geodesic
distance in $S^{d-1}$, in other words the angle formed by $x$, $0$ and $y$.
\begin{lemma}\label{lm:cap}
Let $\Theta \subseteq S^{d-1}$ satisfy $\Theta = - \Theta$, and $\delta>0$.
Suppose also that $K\subseteq\RR^d$ is a convex body
and $D \subseteq \partial K$ is a measurable subset of
the boundary on which almost every point has a normal vector $\xi\in\Theta$.
Let also $\sigma_D$ be the surface measure of $\partial K$ restricted to $D$.
Then
\begin{equation}\label{sigma-ft}
\lim_{t\to\infty} \ft{\sigma_D}(t\eta) = 0,
\end{equation}
uniformly for
$$
\eta \in {\cal N} = \Set{x\in S^{d-1}:~~\dist(x, \Theta)\ge \delta}.
$$
\end{lemma}
{\bf Proof.}
If the set $D$ is contained in a hyperplane
then $\Theta$ may be taken to be a the set $\Set{\pm\theta_0}$.
If, in addition, it is a rectangle, one
gets the validity of \eqref{sigma-ft} by direct calculation of $\ft{\sigma_D}$,
which tends to $0$ except in the normal direction $\theta_0$ and does uniformly so
if one keeps away from the normal direction by any fixed angle.
If $D$ is not a rectangle in its hyperplane one gets \eqref{sigma-ft}
by approximating the set $D$ with a union of finitely
many disjoint rectangles.
Further, if $D$ is polytopal, then it consists
of a finite number of flat pieces $D_i$ and we can represent $\sigma_D$
as the sum of $\sigma_{D_i}$, the surface measure restricted to $D_i$.
Since each of the $\ft{\sigma_{D_i}}(t\eta)$ goes to $0$ uniformly in ${\cal N}$
so does their sum $\ft{\sigma_D}(t\eta)$.

To any measure $\mu \in M(\RR^d)$ and $\eta\in S^{d-1}$ we associate the
projection measure on the line of $\eta$, $\mu^\eta \in M(\RR)$, defined by
\begin{equation}\label{projection-measure}
\mu^\eta(A) = \mu( A\eta + \eta^\perp ),\ \ \ \mbox{($A$ a Borel subset of $\RR$)},
\end{equation}
where $\eta^\perp$ is the hyperplane orthogonal to $\eta$.
By Fubini's theorem we see easily that
\beql{proj-ft}
\ft{\mu^\eta}(t) = \ft{\mu}(t\eta),\ \ \ t\in\RR.
\eeq

We now use the fact that for all $\epsilon,\delta>0$ there
exists a polytopal approximation $P$ of $D$
such that
\begin{itemize}
\item
all normals to $P$ are at distance $\le \delta/2$ from $\Theta$,
\item
for all $\eta\in {\cal N}$ we have that the measures $\sigma_D^\eta$
and $\sigma_P^\eta$ are in fact $L^1$ functions and
\begin{equation}\label{proj-comparison}
\Norm{\sigma_D^\eta - \sigma_P^\eta}_{L^1(\RR)}
	\le\epsilon,\ \ \ (\eta \in {\cal N}).
\end{equation}
\end{itemize}
To prove this one first shows this under the assumption that the
normal vector is a continuous function on $D$ and then uses the well
known theorem (see, for example, \cite[p.\ 23]{morgan}) which says that
we can throw away a part of the surface of arbitrarily small measure
so that the normal is continuous on what remains.

From \eqref{proj-comparison} and from \eqref{proj-ft}
it follows that $\Abs{\ft{\sigma_D}(t\eta) - \ft{\sigma_P}(t\eta)}\le\epsilon$
for all $t\in\RR$, $\eta\in N$.
Since $\epsilon$ is arbitrary Lemma \ref{lm:cap} follows
for the general convex surface piece $D$
from the fact that it holds for polytopal surfaces.
\Qed

\noindent
{\bf Proof of Theorem \ref{th:non-polytopes-afford}.}
(a) Suppose that $K$ is a symmetric polytope with
$N$ non-parallel faces and that $\mu$ is a probability measure on $\partial K$
which is $\epsilon$-good.
We will show that $\epsilon\ge{1\over \sqrt{2} N}$.

The distinct direction vectors $\pm\theta_1,\ldots,\pm\theta_N$
are all the normals that appear on the faces of $K$.
Since all faces are partitioned into those which are normal to $\theta_1$,
normal to $\theta_2$, and so on, it follows that at least one of these pairs
of faces has total $\mu$-measure at least $1/N$,
say the pair of faces normal to $\theta_1$.

The projection measure (see \eqref{projection-measure})
$\mu^{\theta_1}$ has then one or two nonnegative point masses $c_1$ and $c_2$
of total mass at least $1/N$.
But, if $c_i$, $i\in I$, are all the point masses of the finite measure $\mu^{\theta_1}$,
Wiener's Theorem tells us that
$$
\sum_{i\in I} \Abs{c_i}^2 =
	\lim_{T\to\infty}{1\over 2T} \int_{-T}^T\Abs{\ft{\mu^{\theta_1}}(t)}^2~dt.
$$
From \eqref{proj-ft} and the assumption that $\mu$ is $\epsilon$-good
it follows that the right hand side above is $\le \epsilon^2$, from which
it follows that
$$
\epsilon^2 \ge c_1^2 + c_2^2 \ge {1 \over 2 N^2},
$$
since $c_1+c_2 \ge 1/N$.
This is the claimed inequality.

(b) For the remaining part of Theorem \ref{th:non-polytopes-afford}
we must show that whenever $K$ is not a polytope
or is a polytope with at least $N$ non-parallel faces, we can find,
for any $\delta>0$,
a probability measure on $\partial K$ which is $({1\over N}+\delta)$-good.
For this we recall that almost all points $x$ on $\partial K$
have a unique tangent hyperplane $T_x$,  whose outward normal unit vector
we denote by $n(x)$.
This map $x \to n(x)$ is called the Gauss map and, through it, a measure is defined
on $S^{d-1}$, called the area measure $S_K$:
$$
S_K(A) = \sigma_{\partial K}\Set{x\in\partial K:~~n(x) \in A},
\ \ \ \mbox{for any Borel set $A \subseteq S^{d-1}$}.
$$
It is well known that the $0$-symmetric convex body is a polytope
with $N$ pairs of opposite faces if and only if the measure $S_K$ is
a symmetric measure supported on $N$ pairs of opposite points on the sphere $S^{d-1}$.
Therefore the support of the area measure of $K$ contains at least $2N$ points
$\pm\theta_1,\ldots,\pm\theta_N$.
For half of these points, the points $\theta_1,\ldots,\theta_N$, we choose
an open neighborhood $N_i \subset S^{d-1}$, $i=1,\ldots,n$, around
each of them so that all these neighborhoods are disjoint,
and we call $\delta_0>0$ the minimum geodesic distance between any two of them.
Let then $D_i = n^{-1}(N_i) \subseteq \partial K$ be the parts
of the boundary almost all
points of which get mapped into $N_i$ via the Gauss map $n(x)$.
This implies that all points in $D_i$ have a normal in $N_i$,
and the $D_i$ all have positive surface measure.

Define now the probability measure $\mu \in M(\partial K)$
to be an appropriate multiple of surface measure on each $D_i$,
so as to have total mass $1/N$ on each $D_i$.
Call $\mu_i$ the measure $\mu$ resrticted to $D_i$.
From Lemma \ref{lm:cap} it follows that $\ft{\mu_i}(t\eta) \to 0$
as $t\to\infty$
uniformly for all $\eta$ which are distance at least $\delta_0/10$ from $N_i$.
And for all $x$ we have trivially $\Abs{\ft{\mu_i}(x)}\le1/N$.

Let now $\delta>0$ and choose $R>0$ large enought so
that for all $i=1,\ldots,N$
we have $\Abs{\ft{\mu_i}(t\eta)} \le \delta/(N-1)$
if $\Abs{t}>R$ and $\eta$ has distance
more than $\delta_0/10$ from $N_i$.
If now $x\in\RR^d$ is an arbitrary vector with $\Abs{x}>R$
the vector $\eta = {1\over\Abs{x}}x$ can have
distance at most $\delta_0/10$ from at most one
neighborhood $N_i$, say from $N_1$.
If follows that
\begin{eqnarray*}
\Abs{\ft{\mu}(x)} &\le& \Abs{\ft{\mu_1}(x)}+\ldots+\Abs{\ft{\mu_N}(x)}\\
 &\le& {1\over N} + (N-1){\delta \over N-1}\\
 & = & {1\over N} + \delta.
\end{eqnarray*}
\Qed

%%%%%%%%%%%%%%%%%%%%%%%%%%%%%%%%%%%%%%%%%%%%%
% References
%%%%%%%%%%%%%%%%%%%%%%%%%%%%%%%%%%%%%%%%%%%%%
\noindent
\ \\
{\bf Bibliography}
\\

\noindent
{\sc\small
Department of Mathematics, University of Crete, Knossos Ave.,
714 09 Iraklio, Greece.\\
E-mail: {\tt kolount@member.ams.org}
}

\end{document}